\documentclass[11pt]{amsart}
\usepackage{latexsym,amssymb,amsmath,amscd,amscd,amsthm,amsxtra}
\usepackage[dvips]{graphicx}

\usepackage{amsfonts}
\usepackage{amssymb}
\newtheorem{prop}{Proposition}[section]
\newtheorem{thm}[prop]{Theorem}

\newtheorem{corol}[prop]{Corollary}

\newtheorem{rem}[prop]{Remark}

\begin{document}

\title[A direct and simple proof]
 {A direct and simple proof of Jacobi identities for determinants}

\author{Kuihua Yan }

\address{College of Mathematics, Physics and Information
 Science, Zhejiang Normal University, Zhejiang Province, Jinhua City,
 China}

\email{yankh@zjnu.cn}

%\thanks{This work is supported by the
%doctoral scientific research startup foundation of Zhejiang Normal
%University (Grant No. ZC304005089).}

%\thanks{The author was also supported in part by the Research
% Council of Slovenia.}



\date{}

\dedicatory{}

\commby{Kuihua Yan}

%%% ----------------------------------------------------------------------

\begin{abstract}
The Jacobi identities play an important role in constructing the
explicit exact solutions of a broad class of integrable systems in
soliton theory. In the paper, a direct and simple proof of the
Jacobi identities for determinants is presented by employing the
Pl$\ddot{u}$cker relations.
\end{abstract}

%%% ----------------------------------------------------------------------
\maketitle
%%% ----------------------------------------------------------------------

\section{Introduction}
Let $A=\big(a_{_{ij}}\big)_{n\times n}$ be a $n$-order matrix.
Denoted by $M_{ij}\;\big(R_{ij}\equiv (-1)^{i+j}M_{ij}\big)$ the
cofactor (algebraic cofactor) of the matrix entry $a_{_{ij}}$ . The
cofactor (algebraic cofactor) of the minor determinant
$\left|\begin{smallmatrix}
a_{_{ij}}&a_{_{il}}\\
a_{_{kj}}&a_{_{kl}}\\
\end{smallmatrix}\right|$ is denoted by $M_{k\,l}^{i\,j}
\;\big(R_{k\,l}^{i\,j}\equiv (-1)^{k+l+i+j}M_{k\,l}^{i\,j}\big)$,
then the following Jacobi identities \cite{Tak,Hir}
\begin{gather}\label{101}
M_{ii}M_{jj}-M_{ij}M_{ji}=M_{i\,j}^{i\,j}\det A,\quad 1\leq
i,\,j\leq n,
\end{gather}
are valid. Though the Jacobi identities have been proved in
\cite{Tak}, as the author in \cite{Hir} said, looking at the proof
of the general case, it is difficult to understand the Jacobi
identities immediately and the author himself came to understand the
result by checking the formulae using computer algebra, looking for
an alternative proof and applying it to actual problems. Here we
will present a direct proof for the Jacobi identities using the
famous Pl$\ddot{u}$cker relations for determinants.

\section{Pl$\ddot{u}$cker relations}

In this section, let's state the Pl$\ddot{u}$cker relations for
determinants.

\begin{thm}\label{thm2-01}
Let $M$ be a $n\times (n-r)$ matrix and
$a_{_1},\,a_{_2},\,\cdots,\,a_{_{2r}}$ $2r$ $n$-order column
vectors, then
\begin{gather}\label{201}
\sum_{\sigma}
(-1)^{k_{_{1}}+\cdots+k_{_{r}}}\big|\begin{array}{cccc} M &
a_{_{k_{_{1}}}}& \cdots & a_{_{k_{_{r}}}}
\end{array}\big|\cdot
\big|\begin{array}{cccc} M & a_{_{k_{_{r+2}}}} & \cdots &
a_{_{k_{_{2r}}}}
\end{array}\big|=0,
\end{gather}
where $\big\{k_{_{1}},\,k_{_{2}},\,\cdots,\,k_{_{2r}}\big\}$
is a permutation of
$\big\{1,\,2,\,\cdots,\,2r\big\}$ and $\sigma$ is the permutation with
$1\leq k_{_{1}}<\cdots<k_{_{r}}\leq 2r$ and
$k_{_{r+1}}<\cdots<k_{_{2r}}$.
\end{thm}
\vskip -0.4cm
\begin{proof}
Firstly, it is obvious that
\begin{align}
\begin{split}
\left|
\begin{array}{cccccc}
M & 0 & a_{_1} & a_{_2} & \cdots & a_{_{2r}} \\
0 & M & a_{_1} & a_{_2} & \cdots & a_{_{2r}}
\end{array}
\right|&=\left|
\begin{array}{cccccc}
M & -M & 0      & 0      & \cdots & 0 \\
0 & M  & a_{_1} & a_{_2} & \cdots & a_{_{2r}}
\end{array}
\right|\\
&=\left|
\begin{array}{cccccc}
M & 0  & 0      & 0      & \cdots & 0 \\
0 & M  & a_{_1} & a_{_2} & \cdots & a_{_{2r}}
\end{array}
\right|=0.
\end{split}
\end{align}
On the other hand, by the classical Laplace expansion for determinants, it
can be obtained that
\begin{align}\label{203}
\begin{split}
&\left|
\begin{array}{cccccc}
M & 0 & a_{_1} & a_{_2} & \cdots & a_{_{2r}}\\
0 & M & a_{_1} & a_{_2} & \cdots & a_{_{2r}}
\end{array}
\right|\\
=&\sum_{1\leq k_{_{1}}<\cdots<k_{_{r}}\leq 2r
\atop k_{_{r+1}}<\cdots<k_{_{2r}}}
(-1)^{\frac{n(n+1)}{2}+k_{_{1}}+\cdots+k_{_{r}}}\big|\begin{array}{ccccc}
M & a_{_{k_{_{1}}}} & a_{_{k_{_{2}}}} & \cdots & a_{_{k_{_{r}}}}
\end{array}\big|\\
&\hskip 5cm\cdot \big|\begin{array}{ccccc} M & a_{_{k_{_{r+1}}}} &
a_{_{k_{_{r+2}}}} & \cdots & a_{_{k_{_{2r}}}}
\end{array}\big|.
\end{split}
\end{align}
Comparing the above two equations, we have
$$
\sum_{\sigma}
(-1)^{k_{_{1}}+\cdots+k_{_{r}}}\big|\begin{array}{cccc} M &
a_{_{k_{_{1}}}} & \cdots & a_{_{k_{_{r}}}}
\end{array}\big|\cdot
\big|\begin{array}{cccc} M & a_{_{k_{_{r+1}}}}& \cdots &
a_{_{k_{_{2r}}}}
\end{array}\big|=0,
$$
where $\big\{k_{_{1}},\,k_{_{2}},\,\cdots,\,k_{_{2r}}\big\}$
is a permutation of
$\big\{1,\,2,\,\cdots,\,2r\big\}$ and $\sigma$ is the permutation with
$1\leq k_{_{1}}<\cdots<k_{_{r}}\leq 2r$ and
$k_{_{r+1}}<\cdots<k_{_{2r}}$.
 \end{proof}

\begin{corol}\label{corol2-01}
Let $M$ be a $n\times (n-2)$ matrix and $a,\,b,\,c$ and $d$ four
$n$-order column vectors , then
\begin{align}\label{205}
\begin{split}
&\big|\begin{array}{ccc} M & a & b
\end{array}\big|\cdot
\big|\begin{array}{ccc}
M & c & d
\end{array}\big|
-\big|\begin{array}{ccc}
M & a & c
\end{array}\big|\cdot
\big|\begin{array}{ccc}
M & b & d
\end{array}\big|\\
&\hskip 4cm+\big|\begin{array}{ccc} M & a & d
\end{array}\big|\cdot
\big|\begin{array}{ccc}
M & b & c
\end{array}\big|=0.\qquad \qed
\end{split}
\end{align}
\end{corol}
\begin{rem}\label{rem2-01}
The above equation \eqref{205} is the simplest case of the
Pl$\ddot{u}$cker relations \cite{Hir} which plays an important role
in nonlinear dynamics and soliton theory due to Sato's theorem
\cite{Sat,OST}. It shows that many of the differential and
difference equations in mathematical physics are merely disguised
versions of the Pl$\ddot{u}$cker relations. For example, Sato
\cite{Sat,OST} first discovered that the KP equation in bilinear
form
\begin{gather}\label{207}
\big(D_x^4-4D_xD_t+3D_y^2\big)f\cdot f=0
\end{gather}
was nothing but a Pl$\ddot{u}$cker relation, where the Hirota's
bilinear operators $D_t,\,D_x$ and $D_y$ are defined by
\begin{align*}
D_x^mD_t^nf\cdot g&=\big(\partial_{_x}-\partial_{_{x'}}\big)^m
\big(\partial_{_t}-\partial_{_{t'}}\big)^nf(x,\,t)\,g(x',\,t')
\big|_{x'=x,\,t'=t}.
\end{align*}
\end{rem}
\begin{rem}\label{rem2-02}
The other forms of Pl$\ddot{u}$cker relations have been widely
applied to algebraic geometry. For example, a projective embedding
of the Grassmann variety Gr(p, n) can be defined by the quadratic
polynomial equations called ``the Pl$\ddot{u}$cker relations" (to
see \cite{KPRS} and its references).
\end{rem}
\begin{rem}\label{rem2-03}
The Pl$\ddot{u}$cker relations are also exactly relational to the
Maya diagrams and Young diagrams \cite{Hir}. It makes them play a
primary role in Combinatorics, Lie theory and Representation Theory
(to see \cite{Kle,IWa} and their references).
\end{rem}

\section{Jacobi identities for determinants}

\begin{thm}\label{thm3-01}
Let $A=\big(a_{_{ij}}\big)_{n\times n}$ be a $n$-order matrix.
Denoted by $M_{k\,l}^{i\,j}$ the cofactor of the minor determinant
$\left|\begin{smallmatrix}
a_{_{ik}}&a_{_{il}}\\
a_{_{jk}}&a_{_{jl}}\\
\end{smallmatrix}\right|$, then
\begin{gather}\label{301}
M_{k\,l}^{i\,j}M_{s\,r}^{i\,j}-M_{k\,s}^{i\,j}M_{l\,r}^{i\,j}+
M_{k\,r}^{i\,j}M_{l\,s}^{i\,j}=0.
\end{gather}
\end{thm}
\vskip -0.4cm
\begin{proof}
It is no less of generality to consider the case of $i<j$ and
$k<l<s<r$. Denoted by $M$ the $(n-2)\times (n-4)$ submatrix obtained
by eliminating the $i$-th and $j$-th rows and the $k$-th, $l$-th,
$s$-th and $r$-th columns from $A$. The four $(n-2)$-order column
vectors obtained by eliminating the $i$-th and $j$-th components
from the $k$-th, $l$-th, $s$-th and $r$-th column vectors in A are
denoted by $a,\,b,\,c$ and $d$ respectively. Then it is easy to see
that
\begin{align*}
\begin{split}
M_{k\,l}^{i\,j}=(-1)^{k+l}\big|\begin{array}{ccc}
M & a & b
\end{array}\big|;\qquad
&M_{s\,r}^{i\,j}=(-1)^{s+r}\big|\begin{array}{ccc}
M & c & d
\end{array}\big|;\\
M_{k\,s}^{i\,j}=(-1)^{k+s}\big|\begin{array}{ccc}
M & a & c
\end{array}\big|;\qquad
&M_{l\,r}^{i\,j}=(-1)^{l+r}\big|\begin{array}{ccc}
M & b & d
\end{array}\big|;\\
M_{k\,r}^{i\,j}=(-1)^{k+r}\big|\begin{array}{ccc}
M & a & d
\end{array}\big|;\qquad
&M_{l\,s}^{i\,j}=(-1)^{l+s}\big|\begin{array}{ccc}
M & b & c
\end{array}\big|.
\end{split}
\end{align*}
Consequently, by employing the Pl$\ddot{u}$cker relation
\eqref{205}, one has
\begin{align*}
&M_{k\,l}^{i\,j}M_{s\,r}^{i\,j}-M_{k\,s}^{i\,j}M_{l\,r}^{i\,j}+
M_{k\,r}^{i\,j}M_{l\,s}^{i\,j}\\
=&(-1)^{k+l+r+s}\Big(\big|\begin{array}{ccc}
M & a & b
\end{array}\big|\cdot
\big|\begin{array}{ccc}
M & c & d
\end{array}\big|
-\big|\begin{array}{ccc}
M & a & c
\end{array}\big|\cdot
\big|\begin{array}{ccc}
M & b & d
\end{array}\big|\\
&\hskip 5cm+\big|\begin{array}{ccc} M & a & d
\end{array}\big|\cdot
\big|\begin{array}{ccc}
M & b & c
\end{array}\big|\Big)\\
=& 0.
\end{align*}
\end{proof}

\begin{rem}\label{rem3-01}
Note that only the indices are important in the equation
\eqref{301}, so it can also be expressed as
$$
(k,\,l)\,(s,\,r)-(k,\,s)\,(l,\,r)+(k,\,r)\,(l,\,s)=0,
$$
which is nothing but a Pl$\ddot{u}$cker relation.
\end{rem}

\begin{thm}\label{thm3-02}
Let $A=\big(a_{_{ij}}\big)_{n\times n}$ be a $n$-order matrix.
Denoted by $M_{j_{_1}\,\cdots\,j_{_r}}^{i_{_1}\,\cdots\,i_{_r}}$
the cofactor of the minor determinant
$\left|\begin{smallmatrix}
a_{_{i_{_1}j_{_1}}}&a_{_{i_{_1}j_{_2}}}&\cdots&a_{_{i_{_1}j_{_r}}}\\
a_{_{i_{_2}j_{_1}}}&a_{_{i_{_2}j_{_2}}}&\cdots&a_{_{i_{_2}j_{_r}}}\\
\vdots&\vdots&\ddots&\ddots\\
a_{_{i_{_r}j_{_1}}}&a_{_{i_{_r}j_{_2}}}&\cdots&a_{_{i_{_r}j_{_r}}}\\
\end{smallmatrix}\right|$.
Choosing $r$ rows and $2r$ columns from A, the according row and column indices
are denoted by $i_{_1},\,\cdots,\,i_{_r}$
and $j_{_1},\,\cdots,\,j_{_{2r}}$ respectively.
It might as well set that $i_{_1}<\cdots<i_{_r}$ and
$j_{_1}<\cdots<j_{_{2r}}$.
Then
\begin{gather}\label{303}
\sum_{\sigma}
(-1)^{k_{_{1}}+\cdots+k_{_{r}}}
M_{k_{_1}\,\cdots\,k_{_r}}^{i_{_1}\,\cdots\,i_{_r}}
M_{k_{_{r+1}}\,\cdots\,k_{_{2r}}}^{i_{_1}\,\cdots\,i_{_r}}
=0,
\end{gather}
where $\big\{k_{_{1}},\,k_{_{2}},\,\cdots,\,k_{_{2r}}\big\}$
is a permutation of
$\big\{j_{_{1}},\,j_{_{2}},\,\cdots,\,j_{_{2r}}\big\}$
and $\sigma$ is the permutation with
$j_{_{1}}\leq k_{_{1}}<\cdots<k_{_{r}}\leq j_{_{2r}}$ and
$k_{_{r+1}}<\cdots<k_{_{2r}}$.
\end{thm}
\vskip -0.4cm
\begin{proof}
It is completely similar to the proof of the above theorem
\ref{thm3-01} by employing the equation \eqref{201}. So here we omit
it.\end{proof}

\begin{thm}\label{thm3-03}{\bf (Jacobi identity \cite{Tak,Hir})}
Let $A=\big(a_{_{ij}}\big)_{n\times n}$ be a $n$-order matrix.
Denoted by $M_{ij}$ the cofactor of the matrix entry $a_{_{ij}}$ in
A. The cofactor of the minor determinant $\left|\begin{smallmatrix}
a_{_{ij}}&a_{_{il}}\\
a_{_{kj}}&a_{_{kl}}\\
\end{smallmatrix}\right|$ is denoted by $M_{k\,l}^{i\,j}$, then
\begin{gather}\label{305}
M_{ii}M_{jj}-M_{ij}M_{ji}=M_{i\,j}^{i\,j}\det A,\qquad 1\leq
i,\,j\leq n.
\end{gather}
\end{thm}
\vskip -0.4cm
\begin{proof} It is no less generality to condider
the special case with $i=1$ and $j=2$.
Firstly, it is easy to see that

\begin{align*}
M_{1\,1}&=a_{_{22}}M_{1\,2}^{1\,2}
+\sum_{k=3}^{n}(-1)^ka_{_{2k}}M_{2\,k}^{1\,2};\qquad M_{2\,2}
=a_{_{11}}M_{1\,2}^{1\,2}
+\sum_{l=3}^{n}(-1)^la_{_{1l}}M_{1\,l}^{1\,2};\\
M_{1\,2}&=a_{_{21}}M_{1\,2}^{1\,2}
+\sum_{k=3}^{n}(-1)^ka_{_{2k}}M_{1\,k}^{1\,2};\qquad M_{2\,1}
=a_{_{12}}M_{1\,2}^{1\,2}
+\sum_{l=3}^{n}(-1)^la_{_{1l}}M_{2\,l}^{1\,2}.
\end{align*}
Then
\begin{align*}
&M_{11}M_{22}-M_{12}M_{21}\\
=&\big(M_{1\,2}^{1\,2}\big)^2\left|
\begin{smallmatrix}
a_{_{11}}&a_{_{12}}\\
a_{_{21}}&a_{_{22}}
\end{smallmatrix}\right|
+M_{1\,2}^{1\,2}\sum_{l=3}^{n}(-1)^{l+1}M_{1\,l}^{1\,2}\left|
\begin{smallmatrix}
a_{_{12}}&a_{_{1l}}\\
a_{_{22}}&a_{_{2l}}
\end{smallmatrix}\right|\\
&\hskip 1.1cm
+M_{1\,2}^{1\,2}\sum_{k=3}^{n}(-1)^kM_{2\,k}^{1\,2}\left|
\begin{smallmatrix}
a_{_{11}}&a_{_{1k}}\\
a_{_{21}}&a_{_{2k}}
\end{smallmatrix}\right|
+\sum_{k,\,l=3}^{n}(-1)^{k+l}M_{1\,k}^{1\,2}M_{2\,l}^{1\,2}\left|
\begin{smallmatrix}
a_{_{1k}}&a_{_{1l}}\\
a_{_{2k}}&a_{_{2l}}
\end{smallmatrix}\right|.
\end{align*}
On the other hand, by the equation \eqref{301}, it can be obtained
that
\begin{align*}
&\sum_{k,\,l=3}^{n}(-1)^{k+l}M_{1\,k}^{1\,2}M_{2\,l}^{1\,2}\left|
\begin{smallmatrix}
a_{_{1k}}&a_{_{1l}}\\
a_{_{2k}}&a_{_{2l}}
\end{smallmatrix}\right|\\
=&\sum_{3\leq k<l\leq
n}(-1)^{k+l}M_{1\,k}^{1\,2}M_{2\,l}^{1\,2}\left|
\begin{smallmatrix}
a_{_{1k}}&a_{_{1l}}\\
a_{_{2k}}&a_{_{2l}}
\end{smallmatrix}\right|+
\sum_{3\leq l<k\leq n}(-1)^{k+l}M_{1\,k}^{1\,2}M_{2\,l}^{1\,2}\left|
\begin{smallmatrix}
a_{_{1k}}&a_{_{1l}}\\
a_{_{2k}}&a_{_{2l}}
\end{smallmatrix}\right|\\
=&\sum_{3\leq k<l\leq n}(-1)^{k+l}
\big(M_{1\,k}^{1\,2}M_{2\,l}^{1\,2}-M_{1\,l}^{1\,2}M_{2\,k}^{1\,2}\big)
\left|
\begin{smallmatrix}
a_{_{1k}}&a_{_{1l}}\\
a_{_{2k}}&a_{_{2l}}
\end{smallmatrix}\right|\\
=&\sum_{3\leq k<l\leq n}(-1)^{k+l} M_{1\,2}^{1\,2}M_{k\,l}^{1\,2}
\left|
\begin{smallmatrix}
a_{_{1k}}&a_{_{1l}}\\
a_{_{2k}}&a_{_{2l}}
\end{smallmatrix}\right|.
\end{align*}
Therefore, by employing the Laplace expasion for determinants, one has
\begin{align*}
&M_{11}M_{22}-M_{12}M_{21}\\
=&M_{1\,2}^{1\,2}\left[M_{1\,2}^{1\,2}\left|
\begin{smallmatrix}
a_{_{11}}&a_{_{12}}\\
a_{_{21}}&a_{_{22}}
\end{smallmatrix}\right|
+\sum_{l=3}^{n}(-1)^{l+1}M_{1\,l}^{1\,2}\left|
\begin{smallmatrix}
a_{_{12}}&a_{_{1l}}\\
a_{_{22}}&a_{_{2l}}
\end{smallmatrix}\right|\right.\\
&\hskip 1.5cm +\sum_{k=3}^{n}(-1)^kM_{2\,k}^{1\,2}\left|
\begin{smallmatrix}
a_{_{11}}&a_{_{1k}}\\
a_{_{21}}&a_{_{2k}}
\end{smallmatrix}\right|\left.
+\sum_{3\leq k<l\leq n}(-1)^{k+l}
M_{k\,l}^{1\,2}
\left|
\begin{smallmatrix}
a_{_{1k}}&a_{_{1l}}\\
a_{_{2k}}&a_{_{2l}}
\end{smallmatrix}\right|\right]\\
=&M_{1\,2}^{1\,2}\det A.
\end{align*}
\end{proof}

\begin{rem}\label{rem3-02}
The Jacobi identities play an important role in soliton theory. For
example \cite{Hir}, if the solutions to the KP equation or the Toda
lattice equation are expressed as Grammian determinants, their
bilinear equations are nothing but the Jacobi identities.
\end{rem}

\begin{corol}\label{corol3-01}
Let $A=\big(a_{_{ij}}\big)_{2n\times 2n}$ be a $2n$-order
antisymmetric matrix,
then $\det A$ is equal to a perfect square of a polynomial of
its matrix entries $a_{_{ij}}\,(1\leq i,\,j\leq n)$.
\end{corol}
\vskip -0.4cm
\begin{proof}
 Following the above symbols,
it is obvious that $M_{11}=M_{22}=0$ and $ M_{12}=-M_{21}$. Hence
using the Jacobi identities, it can be obtained that
\begin{gather}\label{307}
M_{1\,2}^{1\,2}\det A=\big(M_{12}\big)^2.
\end{gather}

Note that $\det A=a_{_{12}}^2$ when $n=1$, and $M_{1\,2}^{1\,2}$ is
the determinant of a $2(n-1)$-order antisymmetric matrix. Therefore,
it can be deduced that $\det A$ is a perfect square of a polynomial
of $a_{_{ij}}\,(1\leq i,\,j\leq n)$ by the recurrence relation
\eqref{307}.
\end{proof}

\begin{rem}\label{rem3-03}
Recalling the definition of a Pfaffian, the square of a $n$-order
Pfaffian is equal to the determinant of a $2n$-order antisymmetric
matrix. The above corollary in a sense ensures that a Pfaffian is
well defined. On the other hand, by the Pfaffian expression for
determinants, the terms of Jacobi identities \eqref{305} can be
expressed as
\begin{align}\label{309}
\det A&=\big(1,\,2,\,\cdots,\,n,\,
n^*,\,\cdots,\,2^*,\,1^*\big);\\
M_{ii}&=\big(1,\,\cdots,\,\hat{i},\,\cdots,\,n,\,
n^*,\,\cdots,\,\hat{i}^*,\,\cdots,\,1^*\big);\\
M_{ij}&=\big(1,\,\cdots,\,\hat{i},\,\cdots,\,n,\,
n^*,\,\cdots,\,\hat{j}^*,\,\cdots,\,1^*\big);\\
M_{ij}^{ij}&=\big(1,\,\cdots,\,\hat{i},\,\cdots,\,\hat{j},\,\cdots,\,n,\,
n^*,\,\cdots,\,\hat{j}^*,\,\cdots,\,\hat{i}^*,\,\cdots,\,1^*\big),
\end{align}
where the Pfaffian entries are defined by
$(i,\,j)=(i^*,\,j^*)=0,\;(i,\,j^*)=-(j^*,\,i)=a_{_{ij}}$. Another
proof of the Pfaffian version of Jacobi identities can also be found
in \cite{IWa}. Since the Jacobi identities are exactly relational to
the structrue of solutions for soliton equations, the Pfaffian
should also be. In fact, recently, many scholars
\cite{KNH,GNT,HZT,GiNi,LiO,LHZ} have committed themselves to this
domain.

\end{rem}

% ------------------------------------------------------------------------

\subsection*{Acknowledgment}
This work was supported in part by the doctoral scientific research
startup foundation of Zhejiang Normal University under Grant No.
ZC304005089.

% ------------------------------------------------------------------------
%Included for Gather Purpose only:
%input "Xbib.bib"
%\bibliographystyle{amsplain}
%\bibliography{xbib}

\end{document}